\documentclass[12pt]{article}
\usepackage{graphicx}
\usepackage{amsmath}

\newtheorem{theorem}{Theorem}

\newtheorem{corollary}[theorem]{Corollary}

\newtheorem{definition}[theorem]{Definition}

\newtheorem{lemma}[theorem]{Lemma}

\newenvironment{proof}[1][Proof]{\textbf{#1.} }{\ \rule{0.5em}{0.5em}}
\input{tcilatex}

\begin{document}

\title{Topological classification of $\mathbf{Z}_{p}^{m}$ actions on surfaces}
\author{Antonio F. Costa and Sergei M. Natanzon}
\maketitle

\textbf{Abstract.} Let $\widetilde{S}$ be a closed (compact without
boundary) oriented surface with genus $g$, and $G$ be a group isomorphic to $%
\mathbf{Z}_{p}^{m}$, where $p$ is a prime integer. An action of $G$ on $S$
is a pair $(\widetilde{S},f)$, where $f$ is a representation of $G$ in the
group of orientation preserving autohomeomorphisms of $\widetilde{S}$. Two
actions $(\widetilde{S},f)$ and $(\widetilde{S^{\prime }},f^{\prime })$ are
called strongly (resp. weakly) equivalent if there is a homeomorphism$,$ $%
\widetilde{\psi }:\widetilde{S}\rightarrow \widetilde{S}^{\prime },$ sending
the orientation of $\widetilde{S}$ to the orientation of $\widetilde{S}%
^{\prime },$ such that $f^{\prime }(h)=\widetilde{\psi }\circ f(h)\circ 
\widetilde{\psi }^{-1},$ (resp. there is an automorphism $\alpha \in Aut(G)$
such that $f^{\prime }\circ \alpha (h)=\widetilde{\psi }\circ f(h)\circ 
\widetilde{\psi }^{-1}$) for all $h\in G.$ We give the full description of
strong and weak equivalence classes. The main idea of our work is the fact
that a fixed point free action of $\mathbf{Z}_{p}^{m}$ on a surface provides
a bilinear antisymmetric form on $\mathbf{Z}_{p}^{m}.$ For instance, we
prove that the weakly equivalence classes of actions of $G$ on surfaces with
orbit space of genus $g$ are in one to one correspondence with the set of
pairs which consist in a positive integer number $k$, $k\leq m-n,$ $k=(m-n)%
\func{mod}2,$ $g\geq \frac{1}{2}(m-n+k),$ and an orbit of the action of $%
Aut(G)$ on the set of unordered $r$-tuples $[C_{1},...,C_{r}]$ of
non-trivial elements generating a subgroup isomorphic to $\mathbf{Z}_{p}^{n}$
and such that $\sum_{1}^{r}C_{i}=0$. We use this result in describing the
moduli space of complex algebraic curves admitting a group of automorphisms
isomorphic to $\mathbf{Z}_{p}^{m}.$

\bigskip

{\small This research is partially supported by RFBR grant N98-01-00612 NW0
047.008.005 and DGICYT}

\bigskip

\section{Introduction}

The abelian group actions on surfaces constitute a classical subject and it
is studied in the papers [N], [E], [J1], [J2], [Na1], [S], [Z]. In [E],
[J1], [J2], [Z], it is established a connection between the topological
equivalence classes of actions and the second homology of the group that is
acting. But some attempts to use these results for the classification of
abelian actions give wrong results in some cases (compare Remark 4.5 of [E]
with Corollary 12 in our Section 4). The full classification has been found
in the cyclic case by J. Nielsen in [N] and for $\mathbf{Z}_{2}^{m}$ in
[Na1]. In this paper we present a direct way to deal with the topological
classification of $\mathbf{Z}_{p}^{m}$ actions, where $p$ is a prime integer
and we obtain a complete answer ($\mathbf{Z}_{p}^{m}=\mathbf{Z}_{p}\oplus $ $%
\overset{m}{...}\oplus \mathbf{Z}_{p}$ and $\mathbf{Z}_{p}=\mathbf{Z}/p%
\mathbf{Z}$). The main idea of our work is the fact that a fixed point free
action of $\mathbf{Z}_{p}^{m}$ provides a bilinear antisymmetric form on $%
\mathbf{Z}_{p}^{m}.$

Let $\widetilde{S}$ be a closed (compact without boundary) oriented surface
with genus $g=g(\widetilde{S})$, and $G$ be a group isomorphic to $\mathbf{Z}%
_{p}^{m}$, where $p$ is a prime integer. An action of $G$ on $\widetilde{S}$
is a pair $(\widetilde{S},f)$, where $f$ is a representation of $G$ in the
group of orientation preserving autohomeomorphisms of $\widetilde{S}$. Two
actions $(\widetilde{S},f)$ and $(\widetilde{S^{\prime }},f^{\prime })$ are
called strongly equivalent if there is a homeomorphism$,$ $\widetilde{\psi }:%
\widetilde{S}\rightarrow \widetilde{S}^{\prime },$ sending the orientation
of $\widetilde{S}$ to the orientation of $\widetilde{S}^{\prime },$ such
that $f^{\prime }(h)=\widetilde{\psi }\circ f(h)\circ \widetilde{\psi }%
^{-1}, $ for all $h\in G.$ We give the full description of strong
equivalence classes, in particular, in the case of fixed point free actions,
the set of such equivalence classes of actions on surfaces of a given genus
appears to be in bijection with the set of bilinear antisymmetric forms $%
(.,.):G^{\ast }\times G^{\ast }\rightarrow \mathbf{Z}_{p}$, where $G^{\ast }$
is the group of forms of $G$ on $\mathbf{Z}_{p}$ (Theorems 8 and 9)$.$ The
case of actions having elements with fixed points is considered in Theorems
13 and 14.

A motivation for our study is the description of the set of connected
components in the moduli space $M^{p,m}$ of pairs $(C,G)$, where $C$ is a
complex algebraic curve and $G\cong \mathbf{Z}_{p}^{m}$ is a group of
automorphisms of $C$. According to [Na2] the description of connected
components of $M^{p,m}$ is reduced to the description of topological classes
of pairs $(\widetilde{S},K)$ where $K$ is a group of autohomeomorphisms of $%
\widetilde{S}$ and $K$ is isomorphic to $\mathbf{Z}_{p}^{m}.$ We consider
that $(\widetilde{S},K)$ and $(\widetilde{S}^{\prime },K^{\prime })$ are
equivalent if there exist a homeomorphism $\varphi :\widetilde{S}\rightarrow 
\widetilde{S}^{\prime }$ such that $K^{\prime }=\varphi \circ K\circ \varphi
^{-1}.$ These equivalence classes are in one to one correspondence with
classes of weak equivalence (in the terminology of Edmonds [E]). Two actions 
$(\widetilde{S},f)$ and $(\widetilde{S^{\prime }},f^{\prime })$ of a group $%
G\cong \mathbf{Z}_{p}^{m}$ are weakly equivalent if there is a homeomorphism 
$\widetilde{\psi }:\widetilde{S}\rightarrow \widetilde{S}^{\prime }$ and an
automorphism $\alpha \in Aut(G)$ such that $f^{\prime }\circ \alpha (h)=%
\widetilde{\psi }\circ f(h)\circ \widetilde{\psi }^{-1}$ for all $h\in G.$
We prove that the weakly equivalence classes of actions of $G$ on surfaces
with orbit space of genus $g$ are in one to one correspondence with the set
of pairs $(k,Aut(G)[C_{1},...,C_{r}]),$ such that $k\leq m-n,$ $k=(m-n)\func{%
mod}2,$ $g\geq \frac{1}{2}(m-n+k),$ $r\geq n$ and $[C_{1},...,C_{r}]$ is an
unordered r-tuple of non-trivial elements generating $H\cong \mathbf{Z}%
_{p}^{n}$ such that $\sum_{1}^{r}C_{i}=0$ (Theorem 16).

\textbf{Acknowledgment.} The final part of this paper was written during the
authors stay at Max-Planck-Institut f\"{u}r Mathematik in Bonn. We would
like to thank this institution for its support and hospitality. The second
author would like to thank UNED, where this joint work started, for its
hospitality.

\section{Algebraic preliminaries.}

Let us consider the standard lattice $\mathbf{Z}^{2g}=\mathbf{Z}\oplus $ $%
\overset{2g}{...}\oplus \mathbf{Z}$ with the standard basis $%
(e_{i})=((0,...,1^{(i)},...0))$. We define the bilinear antisymmetric form $%
(.,.):\mathbf{Z}^{2g}\times \mathbf{Z}^{2g}\rightarrow \mathbf{Z}$, by $%
(e_{i},e_{j})=\delta _{i+j,2g+1}$ for $i<j$.

We consider also the group $\mathbf{Z}_{p}^{2g}=\mathbf{Z}_{p}\oplus $ $%
\overset{2g}{...}\oplus \mathbf{Z}_{p}$ where $p$ is a prime and $\mathbf{Z}%
_{p}=\{\overline{0},\overline{1},\overline{2},...,\overline{p-1}\}$. Let $%
\varphi :\mathbf{Z}^{2g}\rightarrow \mathbf{Z}_{p}^{2g}$ be the natural
projection defined by $\varphi (e_{i})=\overline{e_{i}}$, where $\overline{%
e_{i}}=(0,...,\overline{1}^{(i)},...0)$. Then we have a bilinear
antisymmetric form $(.,.)_{p}:\mathbf{Z}_{p}^{2g}\times \mathbf{Z}%
_{p}^{2g}\rightarrow \mathbf{Z}_{p}$ defined by $(\overline{e_{i}},\overline{%
e_{j}})=(\varphi (e_{i}),\varphi (e_{j}))_{p}=(e_{i},e_{j})\func{mod}p.$

Let $SL(2g,\mathbf{Z})$ and $SL_{p}(2g,\mathbf{Z}_{p})$ be the subgroups of
the automorphisms groups of $\mathbf{Z}^{2g}$ and $\mathbf{Z}_{p}^{2g}$ that
preserve the bilinear forms $(.,.)$ and $(.,.)_{p}$ respectively. The
natural projection $\varphi :\mathbf{Z}^{2g}\rightarrow \mathbf{Z}_{p}^{2g}$
induces a homomorphism $\varphi _{\ast }:SL(2g,\mathbf{Z})\rightarrow
SL_{p}(2g,\mathbf{Z}_{p})$ such that $\varphi _{\ast }(f)\circ \varphi
=\varphi \circ f$ for all $f\in SL(2g,\mathbf{Z}).$

The following result is elementary:

\begin{theorem}
$\varphi _{\ast }(SL(2g,\mathbf{Z}))=SL_{p}(2g,\mathbf{Z}_{p}).$
\end{theorem}

\begin{proof}
(Sketch). We say that an element of $\mathbf{Z}^{2g}$ is \textit{primitive}
if $e\neq nf$ for all $n\in \mathbf{Z}$ and $f\in \mathbf{Z}^{2g}$. For
every $a_{p}\in \mathbf{Z}_{p}^{2g}$ there is a primitive $a\in \mathbf{Z}%
^{2g}$ such that $\varphi (a)=a_{p}.$

The proof makes use of the two following claims:

1. Assume that $\left\langle \widetilde{a}_{1},\widetilde{a}%
_{2}\right\rangle $ is a subgroup of $\mathbf{Z}_{p}^{2g}$, such that $%
\left\langle \widetilde{a}_{1},\widetilde{a}_{2}\right\rangle \cong \mathbf{Z%
}_{p}^{2}$ and the bilinear form restricted to $\left\langle \widetilde{a}%
_{1},\widetilde{a}_{2}\right\rangle $ is not trivial. Then there are
primitive elements $a_{1},a_{2}\in \mathbf{Z}^{2g}$ such that $%
(a_{1},a_{2})=1$ and $\varphi (a_{1}),\varphi (a_{2})\in \left\langle 
\widetilde{a}_{1},\widetilde{a}_{2}\right\rangle $.

2. Let $a,b\in \mathbf{Z}^{2}$, $\varphi (a)\neq 0,$ $\varphi (b)\neq 0$ and 
$(a,b)=mp$ with $m\in \mathbf{Z}$. Then $\left\langle \varphi
(a)\right\rangle =\left\langle \varphi (b)\right\rangle .$

Using induction and the claims 1 and 2 it is easy to prove:

Let $G$ be a subgroup of $\mathbf{Z}_{p}^{2g}$. Then there is $\Delta
=((a_{i},(i=1,...,r),b_{j},(j=1,...,k\leq s))\subset \mathbf{Z}^{2g}$,
(where $k$ may be 0), such that $\varphi (\Delta )$ generate $G$, $\Delta $
is linear independent and $(a_{i},a_{j})=(b_{i},b_{j})=0$, $%
(a_{i},b_{j})=\delta _{ij}.$

From this fact and induction the Theorem follows.
\end{proof}

We shall also need the following results of linear algebra over finite
fields.

\begin{lemma}
Let $H\cong \mathbf{Z}_{p}^{m}$ and $(.,.):H\times H\rightarrow \mathbf{Z}%
_{p}$ be a bilinear antisymmetric form. Let $\Delta
=(a_{i},(i=1,...,r))\subset H$ be a maximal set of linear independent
elements such that $(a_{i},a_{j})=0$. Then there is a basis of $H,$ $%
((a_{i},(i=1,...,r),b_{j},(j=1,...,s))$, $0\leq s\leq r$, such that $%
(a_{i},a_{j})=(b_{i},b_{j})=0$, $(a_{i},b_{j})=\delta _{ij}.$
\end{lemma}

\begin{proof}
Let us consider all systems $\Delta ^{\prime
}=(a_{i},(i=1,...,r),b_{j},(j=1,...,k))$, $0\leq s\leq r$, such that $%
(a_{i},a_{j})=(b_{i},b_{j})=0$, $(a_{i},b_{j})=\delta _{ij}$. Between them
we choose a system $\Lambda $ with maximal $k.$ Then $\Lambda $ is a basis
with the conditions that we need.
\end{proof}

\begin{theorem}
Let $G,G^{\prime }$ be subgroups of $\mathbf{Z}_{p}^{2g}$ and $\psi
:G\rightarrow G^{\prime }$ be an isomorphism such that $(\psi (a),\psi
(b))_{p}=(a,b)_{p}$ for all $a,b\in G$. Then there is an automorphism $%
\widetilde{\psi }\in SL_{p}(2g,\mathbf{Z}_{p})$ such that $\widetilde{\psi }$
restricted to $G$ is $\psi $.
\end{theorem}

\begin{proof}
It is a consequence of Theorem 1 and Lemma 2.
\end{proof}

\section{Strong classification of fixed point free orientation preserving
actions of $\mathbf{Z}_{p}^{m}$ on surfaces}

Let $\widetilde{S}$ be a closed (compact without boundary) oriented surface
with genus $g$, and $G$ be a group isomorphic to $\mathbf{Z}_{p}^{m}$, where 
$p$ is a prime integer. An action of $G$ on $\widetilde{S}$ is a pair $(%
\widetilde{S},f)$, where $f$ is a monomorphism of $G$ in the group of
orientation preserving autohomeomorphisms of $\widetilde{S}$.

\begin{definition}
(Strong equivalence). Two actions $(\widetilde{S},f)$ and $(\widetilde{%
S^{\prime }},f^{\prime })$ are called strongly equivalent if there is a
homeomorphism$,$ $\widetilde{\psi }:\widetilde{S}\rightarrow \widetilde{S}%
^{\prime },$ sending the orientation of $\widetilde{S}$ to the orientation
of $\widetilde{S}^{\prime }$ and such that $f^{\prime }(h)=\widetilde{\psi }%
\circ f(h)\circ \widetilde{\psi }^{-1},$ for all $h\in G.$
\end{definition}

We are interested in to find all the strong equivalence classes of actions
of $\mathbf{Z}_{p}^{m}$.

We denote by $S=\widetilde{S}/f(G)$ and by $\varphi =\varphi (f):\widetilde{S%
}\rightarrow S$ the natural projection. We shall consider first the case
when $f(h)$ has no fixed points for any $h\in G$, i. e., the action of $(%
\widetilde{S},f)$ is fixed point free. The general case will be considered
in Section 5. Then the projection $\varphi (f):\widetilde{S}\rightarrow S$
is an unbranched covering with deck group of transformations $f(G)$.

Let us consider $\pi _{1}(S)$ as the group of deck transformations of the
universal covering of $S$. Then we have:

\begin{center}
$\omega (\widetilde{S},f):\pi _{1}(S)\rightarrow \pi _{1}(S)/\pi _{1}(%
\widetilde{S})=f(G)\overset{f^{-1}}{\rightarrow }G.$
\end{center}

The resulting epimorphism $\omega (\widetilde{S},f):\pi _{1}(S)\rightarrow
G\cong \mathbf{Z}_{p}^{m}$ is the monodromy epimorphism of the covering $%
\varphi (f):\widetilde{S}\rightarrow S$. The epimorphism $\omega (\widetilde{%
S},f):\pi _{1}(S)\rightarrow G$, induces the epimorphism $\theta _{p}(%
\widetilde{S},f):H_{1}(S,\mathbf{Z}_{p})\rightarrow G$, since $G$ is abelian.

Conversely, given an epimorphism $\theta _{p}:H_{1}(S,\mathbf{Z}%
_{p})\rightarrow G$, there is an action $(\widetilde{S},f)$ such that $%
\theta _{p}=\theta _{p}(\widetilde{S},f).$ To obtain $\widetilde{S}$ is
enough to consider the monodromy $\omega :\pi _{1}(S)\rightarrow H_{1}(S,%
\mathbf{Z}_{p})\overset{\theta _{p}}{\rightarrow }G$ and then $\widetilde{S}%
=U/\ker \omega $, where $U$ is the universal covering of $S$ and the action
of $G$ is given by $G=\pi _{1}(S)/\ker \omega $.

\begin{definition}
Let $S$ and $S^{\prime }$ be two surfaces. Two epimorphisms $\theta :H_{1}(S,%
\mathbf{Z}_{p})\rightarrow G$ and $\theta ^{\prime }:H_{1}(S^{\prime },%
\mathbf{Z}_{p})\rightarrow G$ are called strongly equivalent if there is a
homeomorphism $\psi :S\rightarrow S^{\prime }$ such that induces an
isomorphism $\psi _{p}:H_{1}(S,\mathbf{Z}_{p})\rightarrow H_{1}(S^{\prime },%
\mathbf{Z}_{p})$ such that $\theta =\theta ^{\prime }\circ \psi _{p}$
\end{definition}

\begin{theorem}
(P. A. Smith [S]). Two actions $(\widetilde{S},f)$ and $(\widetilde{S}%
^{\prime },f^{\prime })$ are strongly equivalent if and only if the
epimorphisms $\theta _{p}(\widetilde{S},f)$ and $\theta _{p}(\widetilde{%
S^{\prime }},f^{\prime })$ are strongly equivalent.
\end{theorem}

\begin{definition}
Let $(\widetilde{S},f)$ be an action of $G,$ $S=\widetilde{S}/f(G)$, and $%
\theta =\theta _{p}(\widetilde{S},f):H_{1}(S,\mathbf{Z}_{p})\rightarrow G$
be the epimorphism defined by the action $(\widetilde{S},f)$. Let us
consider the spaces of homomorphisms $G^{\ast }=\{e:G\rightarrow \mathbf{Z}%
_{p}\}$ and $H^{1}(S,\mathbf{Z}_{p})=\{e:H_{1}(S,\mathbf{Z}_{p})\rightarrow 
\mathbf{Z}_{p}\}$. Then $\theta $ generates a monomorphism $\theta ^{\ast
}=\theta ^{\ast }(\widetilde{S},f):G^{\ast }\rightarrow H^{1}(S,\mathbf{Z}%
_{p}).$ The intersection form $(.,.)_{p}=(.,.)_{p}^{S}$ on $H_{1}(S,\mathbf{Z%
}_{p})$ induces an isomorphism $i:H^{1}(S,\mathbf{Z}_{p})\rightarrow H_{1}(S,%
\mathbf{Z}_{p})$ defined by $(a,.)\rightarrow a$ and a form $(.,.)_{(%
\widetilde{S},f)}:G^{\ast }\times G^{\ast }\rightarrow \mathbf{Z}_{p}$ such
that $(a,b)_{(\widetilde{S},f)}=(i\circ \theta ^{\ast }(a),i\circ \theta
^{\ast }(b))_{p}.$
\end{definition}

\begin{theorem}
Two actions $(\widetilde{S},f)$ and $(\widetilde{S}^{\prime },f^{\prime })$
of the group $G$ are strongly equivalent if and only if $\widetilde{S}$ and $%
\widetilde{S}^{\prime }$ have the same genus and $(.,.)_{(\widetilde{S}%
,f)}=(.,.)_{(\widetilde{S^{\prime }},f^{\prime })}.$
\end{theorem}

\begin{proof}
Let $S$ and $S^{\prime }$ denote $\widetilde{S}/f(G)$ and $\widetilde{%
S^{\prime }}/f^{\prime }(G)$ respectively. Assume that $(\widetilde{S},f)$
and $(\widetilde{S}^{\prime },f^{\prime })$ are strongly equivalent, then
according to Theorem 6 there exists a homeomorphism $\psi :S\rightarrow
S^{\prime }$, which induces an isomorphism $\psi _{p}:H_{1}(S,\mathbf{Z}%
_{p})\rightarrow H_{1}(S^{\prime },\mathbf{Z}_{p})$ such that $\theta
=\theta ^{\prime }\circ \psi _{p}$. Since $\psi _{p}$ is induced by a
homeomorphism then $\psi _{p}$ preserves the intersection form and induces
an isomorphism $\psi ^{\ast }:H^{1}(S^{\prime },\mathbf{Z}_{p})\rightarrow
H^{1}(S,\mathbf{Z}_{p})$ such that $(a,b)_{p}^{S^{\prime }}=(\psi ^{\ast
}(a),\psi ^{\ast }(b))_{p}^{S}$ and $\theta ^{\ast }(\widetilde{S},f)=\psi
^{\ast }\circ \theta ^{\ast }(\widetilde{S}^{\prime },f^{\prime })$. Hence,
for $a,b\in G^{\ast }$ we have

$(a,b)_{(\widetilde{S},f)}=(i\circ \theta ^{\ast }(\widetilde{S}%
,f)(a),i\circ \theta ^{\ast }(\widetilde{S},f)(b))_{p}^{S}=$

$=(i\circ \psi ^{\ast }\circ \theta ^{\ast }(\widetilde{S}^{\prime
},f^{\prime })(a),i\circ \psi ^{\ast }\circ \theta ^{\ast }(\widetilde{S}%
^{\prime },f^{\prime })(b))_{p}^{S}=$

$=(i^{\prime }\circ \theta ^{\ast }(\widetilde{S}^{\prime },f^{\prime
})(a),i^{\prime }\circ \theta ^{\ast }(\widetilde{S}^{\prime },f^{\prime
})(b))_{p}^{S^{\prime }}=(a,b)_{(\widetilde{S}^{\prime },f^{\prime })}.$

Assume now $(.,.)_{(\widetilde{S},f)}=(.,.)_{(\widetilde{S^{\prime }}%
,f^{\prime })}.$ We consider some isomorphisms

$Q:H^{1}(S,\mathbf{Z}_{p})\rightarrow (\mathbf{Z}_{p}^{2g},(.,.))$ and $%
Q^{\prime }:H^{1}(S^{\prime },\mathbf{Z}_{p})\rightarrow (\mathbf{Z}%
_{p}^{2g},(.,.))$,

such that $(Q(a),Q(b))_{p}=(a,b)_{p}^{S}$ and $(Q^{\prime }(a^{\prime
}),Q^{\prime }(b^{\prime }))_{p}=(a^{\prime },b^{\prime })_{p}^{S}$, for any 
$a,b\in H^{1}(S,\mathbf{Z}_{p})$ and $a^{\prime },b^{\prime }\in
H^{1}(S^{\prime },\mathbf{Z}_{p})$.

We note $\widetilde{G}=Q\circ \theta ^{\ast }(\widetilde{S},f)(G^{\ast
})\subset \mathbf{Z}_{p}^{2g}$, and $\widetilde{G}^{\prime }=Q^{\prime
}\circ \theta ^{\ast }(\widetilde{S}^{\prime },f^{\prime })(G^{\ast
})\subset \mathbf{Z}_{p}^{2g}.$ Let $\psi :\widetilde{G}\rightarrow 
\widetilde{G}^{\prime }$ be the isomorphism given by $\psi =Q^{\prime }\circ
Q^{-1}$. Then, for every $a,b\in \widetilde{G}$, we have $(\psi (a),\psi
(b))_{p}=(a,b)_{p}.$ From Theorem 3 follows that there is $\widetilde{\psi }%
\in SL_{p}(2g,\mathbf{Z})$ such that $\widetilde{\psi }$ restricted to $%
\widetilde{G}$ is $\psi $. Consider now $\Psi =Q^{-1}\circ \widetilde{\psi }%
\circ Q^{\prime }:H^{1}(S^{\prime },\mathbf{Z}_{p})\rightarrow H^{1}(S,%
\mathbf{Z}_{p}).$ Since $\widetilde{\psi }\in SL_{p}(2g,\mathbf{Z})$ then $%
\Psi $ comes from an isomorphism $\psi _{\ast }:H_{1}(S,\mathbf{Z}%
)\rightarrow H_{1}(S^{\prime },\mathbf{Z})$ sending the intersection form of 
$H_{1}(S,\mathbf{Z})$ to the intersection form of $H_{1}(S^{\prime },\mathbf{%
Z})$ (Theorem 1)$.$ Then by a classical result of H. Burkardt in 1890 (see
[MKS], pg. 178), there is some homeomorphism $\psi :S\rightarrow S^{\prime }$
inducing $\psi _{\ast }$ and $\Psi $, and by construction $\theta ^{\ast }(%
\widetilde{S},f)=\Psi \circ \theta ^{\ast }(\widetilde{S}^{\prime
},f^{\prime }).$ Then by Theorem 6 the actions $(\widetilde{S},f)$ and $(%
\widetilde{S}^{\prime },f^{\prime })$ are strongly equivalent.
\end{proof}

\begin{theorem}
Let $G\cong \mathbf{Z}_{p}^{m}$ and $(.,.):G^{\ast }\times G^{\ast
}\rightarrow \mathbf{Z}_{p}$ be a bilinear antisymmetric form where $k=\dim
\{h\in G^{\ast }:(h,G^{\ast })=0\}$. Then an action $(\widetilde{S},f)$ such
that $(.,.)=(.,.)_{(\widetilde{S},f)}$ and $g=g(\widetilde{S}/f(G))$ exists
if and only if $g\geq \frac{1}{2}(m+k),$ $k=m\func{mod}2,$ $k\leq m$.
\end{theorem}

\begin{proof}
First we construct the action from the form and the numerical conditions.
Applying Lemma 2 we have a basis of $G^{\ast }$, $(a_{i}^{\ast
},(i=1,...,r),b_{j}^{\ast },(j=1,...,s))$ $0\leq s\leq r$, such that $%
(a_{i}^{\ast },a_{j}^{\ast })=(b_{i}^{\ast },b_{j}^{\ast })=0$, $%
(a_{i}^{\ast },b_{j}^{\ast })=\delta _{ij}$ and $s-r=k$. Let $%
(a_{i},(i=1,...,r),b_{j},(j=1,...,s))$ be the dual basis of the above one.
Now consider a surface $S$ of genus $g$ and a basis of $H_{1}(S,\mathbf{Z}%
_{p}),$ $(\alpha _{i},(i=1,...,g),\beta _{i},(i=1,...,g))$. Then we
construct the epimorphism $\theta :H_{1}(S,\mathbf{Z}_{p})\rightarrow
G^{\ast }$, defined by $\theta (\alpha _{i})=a_{i},$ if $i\leq r,$ $\theta
(\alpha _{i})=0,$ if $i>r$ and $\theta (\beta _{i})=b_{i},$ if $i\leq s,$ $%
\theta (\beta _{i})=0,$ if $i>s.$ Then the epimorphism $\theta $ defines a
regular covering $\widetilde{S}\rightarrow S$ with automorphism group $G$
and the action of $G$ on $\widetilde{S}$ satisfies $(.,.)_{(\widetilde{S}%
,f)}=(.,.)$.

Conversely if there is an action $(\widetilde{S},f)$ such that $%
(.,.)=(.,.)_{(\widetilde{S},f)}$ it is obvious that $g=g(\widetilde{S}%
/f(G))\geq \frac{1}{2}(m+k),$ $k=m\func{mod}2,$ $k\leq m$.
\end{proof}

\section{Weak classification of fixed point free orientation preserving
actions of $\mathbf{Z}_{p}^{m}$ on surfaces}

\begin{definition}
(Weak equivalence) Let $(\widetilde{S},f)$ and $(\widetilde{S^{\prime }}%
,f^{\prime })$ be two actions of a group $G\cong \mathbf{Z}_{p}^{m}.$ We
shall say that $(\widetilde{S},f)$ and $(\widetilde{S^{\prime }},f^{\prime
}) $ are weakly equivalent if there is a homeomorphism $\widetilde{\psi }:%
\widetilde{S}\rightarrow \widetilde{S}^{\prime }$ and an automorphism $%
\alpha \in Aut(G)$ such that $f^{\prime }\circ \alpha (h)=\widetilde{\psi }%
\circ f(h)\circ \widetilde{\psi }^{-1}$, $h\in G.$
\end{definition}

The next Theorem solves the problem of weak classification of actions of $%
\mathbf{Z}_{p}^{m}$ on surfaces:

\begin{theorem}
Let $(\widetilde{S},f)$ and $(\widetilde{S^{\prime }},f^{\prime })$ be two
actions of a group $G\cong \mathbf{Z}_{p}^{m}.$ Let $(.,.)_{(\widetilde{S}%
,f)}$ and $(.,.)_{(\widetilde{S^{\prime }},f^{\prime })}$ be the
antisymmetric forms induced by the two actions, $k(\widetilde{S},f)=\dim
\{h\in G^{\ast }:(h,G^{\ast })_{(\widetilde{S},f)}=0\}$ and $k(\widetilde{S}%
^{\prime },f^{\prime })=\dim \{h\in G^{\ast }:(h,G^{\ast })_{(\widetilde{S}%
^{\prime },f^{\prime })}=0\}$. Then the actions $(\widetilde{S},f)$ and $(%
\widetilde{S^{\prime }},f^{\prime })$ are weakly equivalent if and only if $%
g(\widetilde{S}/f(G))=g(\widetilde{S}^{\prime }/f^{\prime }(G))$ and $k(%
\widetilde{S},f)=k(\widetilde{S}^{\prime },f^{\prime }).$
\end{theorem}

\begin{proof}
Let us call $S=\widetilde{S}/f(G)$ and $S^{\prime }=\widetilde{S}^{\prime
}/f^{\prime }(G),$ $g=$ $g(S)$ and $g^{\prime }=g(S^{\prime })$. Let $\theta
^{\ast }(\widetilde{S},f)$ and $\theta ^{\ast }(\widetilde{S}^{\prime
},f^{\prime })$ be the epimorphisms defined by the two actions, $\widetilde{G%
}$ be the image of $G^{\ast }$ in $H_{1}(S,\mathbf{Z}_{p})$ by $\theta
^{\ast }(\widetilde{S},f)$ and $\widetilde{G}^{\prime }$ be the image of $%
G^{\ast }$ in $H_{1}(S^{\prime },\mathbf{Z}_{p})$ by $\theta ^{\ast }(%
\widetilde{S}^{\prime },f^{\prime })$.

Assume that $g(S)=g(S^{\prime })$ and $k(\widetilde{S},f)=k(\widetilde{S}%
^{\prime },f^{\prime })$. Since $k(\widetilde{S},f)=k(\widetilde{S}^{\prime
},f^{\prime })$ then there exists an isomorphism $\psi :\widetilde{G}%
^{\prime }\rightarrow \widetilde{G}$ such that

$(\psi (a),\psi (b))_{(\widetilde{S}^{\prime },f^{\prime })}=(a,b)_{(%
\widetilde{S},f)}.$

Then, using Theorem 3, and that $g(S)=g(S^{\prime })$, there is an
isomorphism $\widetilde{\psi }:H^{1}(S^{\prime },\mathbf{Z}_{p})\rightarrow
H^{1}(S,\mathbf{Z}_{p})$ giving by restriction $\psi $ and sending the
intersection form of $H^{1}(S^{\prime },\mathbf{Z}_{p})$ to the intersection
form of $H^{1}(S,\mathbf{Z}_{p})$. By [MKS, pag 178], there exists a
homeomorphism $\varphi :S\rightarrow S^{\prime }$ inducing $\widetilde{\psi }
$ on cohomology. Then by Theorem 6, the actions $(\widetilde{S},f)$ and $(%
\widetilde{S},\varphi ^{-1}\circ f^{\prime }\circ \varphi )$ are strongly
equivalent. The isomorphism $\psi $, defines an automorphism of $G$ giving
the weak equivalence between $(\widetilde{S},f)$ and $(\widetilde{S^{\prime }%
},f^{\prime }).$
\end{proof}

\begin{corollary}
The weak equivalence classes of $\mathbf{Z}_{p}^{m}$ actions are in
bijection with the set of pairs of positive integer numbers $(k,g)$ such
that $k\leq m,$ $k=m\func{mod}2$ and $g\geq \frac{1}{2}(m+k).$
\end{corollary}

\section{Classification of orientation preserving actions of $\mathbf{Z}%
_{p}^{m}$ with elements having fixed points.}

Let $G$ be a group isomorphic to $\mathbf{Z}_{p}^{m}$ and $(\widetilde{S},f)$
be an action of $G$ on an oriented closed surface $\widetilde{S}$ $.$ We
shall call $G_{fix}$ to the subgroup of $G$ generated by the elements of $%
f(G)$ having fixed points.

The projection $\varphi =\varphi (f):\widetilde{S}\rightarrow S=\widetilde{S}%
/f(G)$ is a covering branched on a finite set of points $\mathcal{B}%
=\{b_{1},...,b_{r}\}$. The covering $\varphi $ is now determined by an
epimorphism $\theta _{p}(\widetilde{S},f):H_{1}(S-\mathcal{B},\mathbf{Z}%
_{p})\rightarrow G$.

We shall call $X_{i}$, $i=1,...,r$, the element of $H_{1}(S-\mathcal{B},%
\mathbf{Z}_{p})$ represented by the boundary of a small disc in $S$ around
the branched point $b_{i},$ and with the orientation given by the
orientation of $S$. Then the set $\{\theta _{p}(\widetilde{S},f)(X_{i})\}$
is a topological invariant for the action $(\widetilde{S},f)$ and $%
\left\langle \theta _{p}(\widetilde{S},f)(X_{i}),i=1,...,r\right\rangle
=G_{fix}.$ Then we have an epimorphism $\vartheta :H_{1}(S,\mathbf{Z}%
_{p})\rightarrow G_{free}$, defined by

\begin{center}
$H_{1}(S,\mathbf{Z}_{p})\rightarrow H_{1}(S-\mathcal{B},\mathbf{Z}%
_{p})/\left\langle X_{i},i=1,...,r\right\rangle \rightarrow
G/G_{fix}=G_{free}$.
\end{center}

In fact the epimorphism $\vartheta $ is the epimorphism defined by the fixed
point free action defined by the unbranched covering $\widetilde{S}%
/f(G_{fix})\rightarrow S$. If $G_{free}=G/G_{fix}$ then $\vartheta $
defines, as in Section 2, a bilinear form $(.,.)_{(\widetilde{S}%
,f)}:G_{free}^{\ast }\times G_{free}^{\ast }\rightarrow \mathbf{Z}_{p}.$

\begin{theorem}
Two actions $(\widetilde{S},f)$ and $(\widetilde{S^{\prime }},f^{\prime })$
of the group $G\cong \mathbf{Z}_{p}^{m}$ are strongly equivalent if and only
if:

1. $\widetilde{S}$ and $\widetilde{S^{\prime }}$ have the same genus, the
number of branched points $r=\#\mathcal{B}$, of the covering $\widetilde{S}%
\rightarrow S=\widetilde{S}/f(G),$ is the same than the number of branched
points $r^{\prime }=\#\mathcal{B}^{\prime }$, of the covering $\widetilde{S}%
^{\prime }\rightarrow S^{\prime }=\widetilde{S}^{\prime }/f^{\prime }(G).$

2. $[\theta _{p}(\widetilde{S},f)(X_{1}),\theta _{p}(\widetilde{S}%
,f)(X_{2}),...,\theta _{p}(\widetilde{S},f)(X_{r})]=$

$=[\theta _{p}(\widetilde{S},f)(X_{1}),\theta _{p}(\widetilde{S}%
,f)(X_{2}),...,\theta _{p}(\widetilde{S},f)(X_{r})]$, where $[.,...,.]$
means unordered r-tuple. As a consequence $G_{free}^{f}=G_{free}^{f^{\prime
}}=G_{free}.$

3. The intersection form on $G_{free},$ induced by $f$ and $f^{\prime }$ is
the same, $(.,.)_{(\widetilde{S},f)}=(.,.)_{(\widetilde{S}^{\prime
},f^{\prime })}.$
\end{theorem}

\begin{proof}
Using Dehn twists along curves around the branched points (see [C], pg. 151,
move (6)) it is possible to obtain a basis $%
(A_{i},(i=1,...,g),B_{i},(i=1,...,g),X_{i},(i=1,...,r)),$ of $H_{1}(S-%
\mathcal{B},\mathbf{Z}_{p})$ such that

$\theta _{p}(\widetilde{S},f)(A_{i})\in G_{free},\theta _{p}(\widetilde{S}%
,f)(B_{i})\in G_{free},$ $i=1,...,g,$

and $(A_{i},A_{j})=0,$ $(B_{i},B_{j})=0,$ $(A_{i},B_{j})=\delta _{ij}.$

In the same way, we can construct a basis $(A_{i}^{\prime
},(i=1,...,g),B_{i}^{\prime },(i=1,...,g),X_{i}^{\prime },(i=1,...,r)),$ of $%
H_{1}(S^{\prime }-\mathcal{B}^{\prime },\mathbf{Z}_{p})$ such that

$\theta _{p}(\widetilde{S}^{\prime },f^{\prime })(A_{i}^{\prime })\in
G_{free},\theta _{p}(\widetilde{S}^{\prime },f^{\prime })(B_{i}^{\prime
})\in G_{free},$ $i=1,...,g,$

and $(A_{i}^{\prime },A_{j}^{\prime })=0,$ $(B_{i}^{\prime },B_{j}^{\prime
})=0,$ $(A_{i}^{\prime },B_{j}^{\prime })=\delta _{ij},$

remark that by condition 1, $g=g^{\prime }$ and $r=r^{\prime }.$

By condition 3 and Theorem 8, then the fixed point free action of $G_{free}$
on $\widetilde{S}/f(G_{fix})$ given by $f$ and the fixed point free action
of $G_{free}$ on $\widetilde{S}^{\prime }/f^{\prime }(G_{fix})$ given by $%
f^{\prime }$ are strongly equivalent. Then there exists a homeomorphism,
preserving the orientations, $\varphi :S\rightarrow S^{\prime }$ inducing on
homology an isomorphism $\psi :H_{1}(S,\mathbf{Z}_{p})\rightarrow
H_{1}(S^{\prime },\mathbf{Z}_{p})$ and by the proof of the Theorem 8 we can
construct $\varphi $ such that $\psi (A_{i})=A_{i}^{\prime }$, and $\psi
(B_{i})=B_{i}^{\prime }$. We now consider a disc $D$ on $S$ containing $%
\mathcal{B}$ and a disc $D^{\prime }$ on $S^{\prime }$ containing $\mathcal{B%
}^{\prime }$. Then we can modify $\varphi $ by composing with an isotopy in $%
S^{\prime }$ in order that $\varphi (D)=D^{\prime }$, and $\varphi
(b_{i})=b_{\sigma (i)}$ where $\sigma $ is a permuntation of $\{1,...,r\}$
such that

$(\theta _{p}(\widetilde{S},f)(X_{1}),\theta _{p}(\widetilde{S}%
,f)(X_{2}),...,\theta _{p}(\widetilde{S},f)(X_{r}))=$

$=(\theta _{p}(\widetilde{S},f)(X_{\sigma (1)}),\theta _{p}(\widetilde{S}%
,f)(X_{\sigma (2)}),...,\theta _{p}(\widetilde{S},f)(X_{\sigma (r)})).$

Now $\varphi $ defines an isomorphim $\widetilde{\psi }:H_{1}(S-\mathcal{B},%
\mathbf{Z}_{p})\rightarrow H_{1}(S^{\prime }-\mathcal{B}^{\prime },\mathbf{Z}%
_{p})$ such that $\widetilde{\psi }(A_{i})=A_{i}^{\prime }$, $\widetilde{%
\psi }(B_{i})=B_{i}^{\prime },$ $\widetilde{\psi }(X_{i})=X_{i}^{\prime },$
then $\theta _{p}(\widetilde{S},f)=\theta _{p}(\widetilde{S}^{\prime
},f^{\prime })\circ \varphi $. Hence the actions $(\widetilde{S},f)$ and $(%
\widetilde{S}^{\prime },f^{\prime })$ are strongly equivalent.
\end{proof}

As a consequence of Theorem 13 and Theorem 9 we have:

\begin{theorem}
Let $G\cong \mathbf{Z}_{p}^{m},$ and $H\cong \mathbf{Z}_{p}^{n}$ be a
subgroup of $G.$ Assume that $[C_{1},...,C_{r}]$, $r\geq n$, be an unordered
element of $(H-\{0\})^{r}$, where $\{C_{1},...,C_{r}\}$ generates $H$ and $%
\sum_{1}^{r}C_{i}=0.$ Let $(.,.)$ be a bilinear antisymmetric form on $G/H$
with $k=\dim \{h\in (G/H)^{\ast }:(h,(G/H)^{\ast })=0\}$. Then for $g\geq 
\frac{1}{2}(m-n+k)$ and only for all such $g$ there is an action $(%
\widetilde{S},f)$ with $g=$ $g(\widetilde{S}/f(G)),$ $(.,.)=(.,.)_{(%
\widetilde{S},f)},$ where $(.,.)_{(\widetilde{S},f)}$ is the bilinear form
induced by the fixed point free action on $G/H,$ the elements acting with
fixed points generate $H$ and $[\theta _{p}(\widetilde{S},f)(X_{1}),\theta
_{p}(\widetilde{S},f)(X_{2}),...,\theta _{p}(\widetilde{S}%
,f)(X_{r})]=[C_{1},C_{2},...,C_{r}]$.
\end{theorem}

\textbf{Remark}. The unordered elements of $H^{r}$, are in one to one
correspondence with the functions $F:H\rightarrow (\mathbf{Z}_{+})^{p-1}$.
From $[C_{1},C_{2},...,C_{r}]$ we define $F(h)=(k_{1},k_{2},...,k_{p-1})$ if
the element $h^{i}$ appears $k_{i}$ times in $[C_{1},C_{2},...,C_{r}]$. The
function $F$ gives the topological type of the action of $H$.

\begin{theorem}
Two actions $(\widetilde{S},f)$ and $(\widetilde{S^{\prime }},f^{\prime })$
of the group $G\cong \mathbf{Z}_{p}^{m}$ are weakly equivalent if and only
if:

1. $\widetilde{S}$ and $\widetilde{S^{\prime }}$ have the same genus, the
number of branched points $r=\#\mathcal{B}$ of the covering $\widetilde{S}%
\rightarrow S=\widetilde{S}/f(G),$ is the same than the number of branched
points $r^{\prime }=\#\mathcal{B}^{\prime }$ of the covering $\widetilde{S}%
^{\prime }\rightarrow S^{\prime }=\widetilde{S}^{\prime }/f^{\prime }(G).$

2. $(\theta _{p}(\widetilde{S},f)(X_{1}),\theta _{p}(\widetilde{S}%
,f)(X_{2}),...,\theta _{p}(\widetilde{S},f)(X_{r}))=$

$=(\gamma \circ \theta _{p}(\widetilde{S},f)(X_{\sigma (1)}),\gamma \circ
\theta _{p}(\widetilde{S},f)(X_{\sigma (2)}),...,\gamma \circ \theta _{p}(%
\widetilde{S},f)(X_{\sigma (r)}))$, where $\sigma $ is a permutation of $%
\{1,...,r\}$ and $\gamma $ is an automorphism of $G$.

3. $\dim \{h\in G_{free}^{\ast }:(h,G_{free}^{\ast })_{(\widetilde{S}%
,f)}=0\}=\dim \{h\in G_{free}^{\ast }:(h,G_{free}^{\ast })_{(\widetilde{S}%
^{\prime },f^{\prime })}=0\}.$
\end{theorem}

\begin{proof}
It is similar to the proof of Theorem 13 but using Theorem 8.
\end{proof}

As a consequence of Corollary 12 we have:

\begin{theorem}
Let $G\cong \mathbf{Z}_{p}^{m}.$ The weak equivalence classes of actions of $%
G$ are in bijection with the set of triples

$(k,g,Aut(G)[C_{1},...,C_{r}]),$

such that $k\leq m-n,$ $k=(m-n)\func{mod}2,$ $g\geq \frac{1}{2}(m-n+k),$ $%
r\geq n$ and $[C_{1},...,C_{r}]$ is an unordered r-tuple of non-trivial
elements of $H^{r}$ such that $\{C_{1},...,C_{r}\}$ generates a group
isomorphic to $\mathbf{Z}_{p}^{n}$ and $\sum_{1}^{r}C_{i}=0.$
\end{theorem}

Let $M^{p,m}$ be the space of pairs $(\widetilde{R},G),$ where $\widetilde{R}
$ is a Riemann surface and $G$ is a group of automorphims of $\widetilde{R}$%
. The covering $\widetilde{R}\rightarrow \widetilde{R}/G$ defines a
projection $p:M^{p,m}\rightarrow M,$ where $M$ is the moduli space of
Riemann surfaces$.$ The projection $p:M^{p,m}\rightarrow M$, gives a
topology on $M^{p,m},$ the weakest topology where $p$ is continuous.

From Theorem 16 and [Na2], Section 6, we have:

\textbf{Consequence} \textit{There exists a one to one correspondence
between the connected components of }$M^{p,m}$\textit{\ with such topology
and triples }

\begin{center}
$(k,g,Aut(G)[C_{1},...,C_{r}])$
\end{center}

\textit{described in Theorem 16. Each connected component of }$M^{p,m}$%
\textit{\ is homeomorphic to the quotient }$R^{n}/Mod$\textit{\ of a vector
space }$R^{n}$\textit{\ by the discontinuous action of a group }$Mod.$

\QTP{Body Math}
$\bigskip $

\begin{center}
{\LARGE References}
\end{center}

[C] A. F. Costa, Classification of the orientation reversing homeomorphisms
of finite order of surfaces, \textit{Topology and its Applications}, \textbf{%
62} (1995) 145-162.

[E] A. L. Edmonds, Surface symmetry I, \textit{Michigan Math. J.}, \textbf{29%
}, (1982), 171-183.

[J1] S. A. Jassim, Finite abelian coverings, \textit{Glasgow Math. J}., 
\textbf{25}, (1984), 207-218.

[J2] S. A. Jassim, Finite abelian actions on surfaces, \textit{Glasgow Math.
J}., \textbf{35}, (1993), 225-234.

[MKS] W. Magnus, A. Karras and D. Solitar, Combinatorial group theory,
Dover, New York, 1976.

[N] J. Nielsen, Die Struktur periodischer Transformationen von Fl\"{a}chen, 
\textit{Danske Vid. Selsk, Mat.-Fys. Medd}. \textbf{15} (1937), 1-77.

[Na1] S. M. Natanzon, Moduli spaces of complex algebraic curves with
isomorphic to $(\mathbf{Z}/2\mathbf{Z})^{m}$ groups of automorphisms, 
\textit{Differential Geometry and its Applications} \textbf{5}, (1995) 1-11.

[Na2] S. M. Natanzon, Moduli of Riemann surfaces, Hurwitz-type spaces, and
their superanalogues, \textit{Russian Math. Surveys}, \textbf{54}:1 (1999)
61-117.

[S] P. A. Smith, Abelian actions of periodic maps on compact surfaces, 
\textit{Michigan Math. J.} \textbf{14} (1967), 257-275.

[Z] B. Zimmermann, Finite abelian group actions on surfaces, \textit{%
Yocohama Mathematical Journal,} \textbf{38}, (1990) 13-21.

\bigskip

{\small Antonio F. Costa}

{\small Departamento de Matem\'{a}ticas Fundamentales,UNED}

28040-{\small Madrid, Spain.}

{\small e-mail:\ acosta@mat.uned.es}

\bigskip

{\small Sergei Natanzon}

{\small Moscow State University and Independent Moscow University,}

{\small Moscow, Russia.}

{\small e-mail: natanzon@mccme.ru}

\end{document}